\documentclass[twocolumn]{article} 
\usepackage{amsfonts,amssymb,amsmath,amscd}
\usepackage{graphicx}

\RequirePackage[alphabetic,msc-links]{amsrefs}

\RequirePackage{hyperref} 
\hypersetup{pdfborder = 0 0 0}

\title{
	Remembering Mark Sapir
}

\author{
	Jean-Camille Birget, Gili Golan, Alexander Olshanskii, and Mikhail Volkov 
	\vspace{1.1em} \\
	\textit{with contributions by} \\
	\vspace{0.4em}\\
	\parbox{\textwidth}{
		\centering
		Victor Guba, Olga Kharlampovich, Simon M. Goberstein, Stuart W. Margolis,  \\
		Lev Shneerson, John Meakin, Ilya Kapovich, Eugene Plotkin, 
		Michael Mihalik,\\ Samuel Corson, and Arman Darbinyan
	}
}
\date{}
\begin{document}
		
	\maketitle
	
		\begin{figure}[h]
		\centering
		\includegraphics[width=0.47\textwidth]{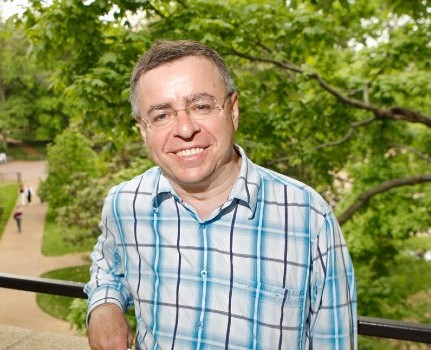}
		\caption{Mark Sapir in Vanderbilt University, 2010.}
		\label{fig:Mark}
	\end{figure}

	Mark Valentinovich Sapir, who passed away in October 2022 at age 65, was born in Sverdlovsk (now Yekaterinburg) in February 1957.  Mark was a brilliant mathematician, whose most significant research contributions were in the areas of geometric group theory, semigroup theory and combinatorial algebra. He was an insightful researcher who achieved remarkable results at the junction of different branches of mathematics.

	Mark entered the Faculty of Mathematics and Mechanics (Matmech) of Ural State University in 1973. As a child prodigy, Mark jumped over two grades of high school so he was much younger than his fellow students. Not only this, but his mathematical talent clearly distinguished Mark among them. Émil Abramovich Golubov 
	who taught first-year students  analytic geometry, was the first to offer Mark research problems. The problems came from the area of Golubov’s immediate interests and were related to various properties of residually finite semigroups.	The collaboration of Sapir and Golubov led to several important results, including their description of varieties consisting entirely of residually finite semigroups.

	After graduating from Matmech, Mark Sapir became a graduate student at the Chair of Algebra and Geometry under the supervision of Lev Naumovich Shevrin, head of the Chair and longtime leader of the research seminar ‘Algebraic Systems’. That seminar played a great role in Sapir’s growth as a mathematician. It is there that Mark presented all his results of the USSR period, giving 65 talks to a well-prepared, demanding, yet friendly audience and gradually developing into a skillful lecturer.
	
	Shevrin suggested that Mark take up the study of semigroup quasivarieties, an area that had not yet been systematically studied at that time. Mark quickly produced an impressive series of results incorporated in his candidate (PhD) dissertation 
	that he defended in 1983. 
	The dissertation made breakthroughs in studying quasi-identity bases and quasivarietal lattices.
	
	During his studies, Mark worked as a teaching assistant in the Department of Physics at Ural State University. It was there that he met Olga, a remarkable student in his class, and they married in 1985.
	
	In 1984, Mark took on the role of Deputy Director at the Computer Science Laboratory of the Sverdlovsk Pedagogical Institute. He co-authored the textbook \textit{Informatics}, designed to teach high school students simple programming, which was printed in over 2.5 million copies. Mark developed the educational software to accompany the textbook, featuring the characters Roo, a kangaroo, and Robby, a snake, who could understand simple  languages. As part of his duties at the pedagogical institute, Mark taught classes to high school students all over the city. Additionally, on his own volition, he organized an afternoon school for high school students, where he used the Roo and Robby software.
	
	During this period, Mark's research was done at night, often in the basement of the pedagogical institute. Despite these challenging conditions, he produced ground-breaking results, including his solutions to Burnside-type problems for semigroups (see below) and a characterization of inherently non-finitely based finite semigroups.
	
	In 1986, Mark and Olga’s  daughter, Jenya, was born. The following year, the Soviet Union began allowing citizens to travel abroad, and Mark took this opportunity to attend international conferences. During one of these trips, he reconnected with Simon Goberstein, who facilitated a visiting position for Mark at California State University, Chico, for eight months in 1991. 
	
	While in Chico, Mark worked on his second dissertation (the doctorate), which was required for securing a permanent research position in the Soviet Union. 
	He initially planned to return to Russia to defend this thesis, but before his stay in Chico ended, the Soviet Union collapsed. With academic prospects in Russia uncertain, particularly because of the quota system limiting opportunities for Jewish scholars, Mark accepted an offer, arranged by John Meakin, to join the University of Nebraska. 
	
	In Nebraska, Mark continued working on semigroups, while beginning the transition to geometric group theory. He collaborated with many of his colleagues including John Meakin, Jean-Camille Birget and Stuart Margolis as well as with Victor Guba, Eliyahu Rips, Alexander Olshanskii and Stanislav Kublanovsky who visited Lincoln. 
	
	In 1993, Sapir and Meakin introduced the concept of diagram groups, which are groups associated with finite semigroup presentations. Mark further developed these ideas in his joint work with Guba. Diagram groups opened up a new framework for studying Thompson's group $F$, a key example of such a group, which remained an important area of Mark’s research until his passing.
	
	During his time in Nebraska,  
	Mark started working on his two most acclaimed papers, which were published in the \textit{Annals of Mathematics} in 2002.	
	
	The first paper, co-authored with Jean-Camille Birget and Eliyahu Rips,
	provided an essentially complete description of all possible sufficiently large growth types of Dehn functions of finitely presented groups. The second paper, joint with Birget, Rips and Olshanskii, proved
	that the word problem of a finitely generated group $G$ is in NP if and only if it can be embedded in a finitely presented group $H$ with polynomial Dehn function.
	
	A key ingredient in these papers was Sapir's invention of S-machines (see below) which is a crucial bridging tool between nondeterministic Turing machines and finite presentations of groups.
	
	In 1997, Mark moved to Vanderbilt University in Nashville, TN,  at the invitation of Ralph Mckenzie. There, he was named a Centennial Professor in 2001. 		By then, Mark focused  primarily on geometric group theory and led a vibrant seminar in the subject. Together with Olshanskii (who moved to Vanderbilt shortly after Sapir), he constructed the first known finitely presented counterexamples to the von Neumann conjecture. 
	 Another notable achievement is his seminal work, primarily with Cornelia Dru\textcommabelow{t}u,  on developing the asymptotic cone approach to the study of relatively hyperbolic groups and his work on lacunary hyperbolic groups, joint with Olshanskii and Osin.

		\begin{figure}[h]
		\centering
		\includegraphics[width=0.47\textwidth]{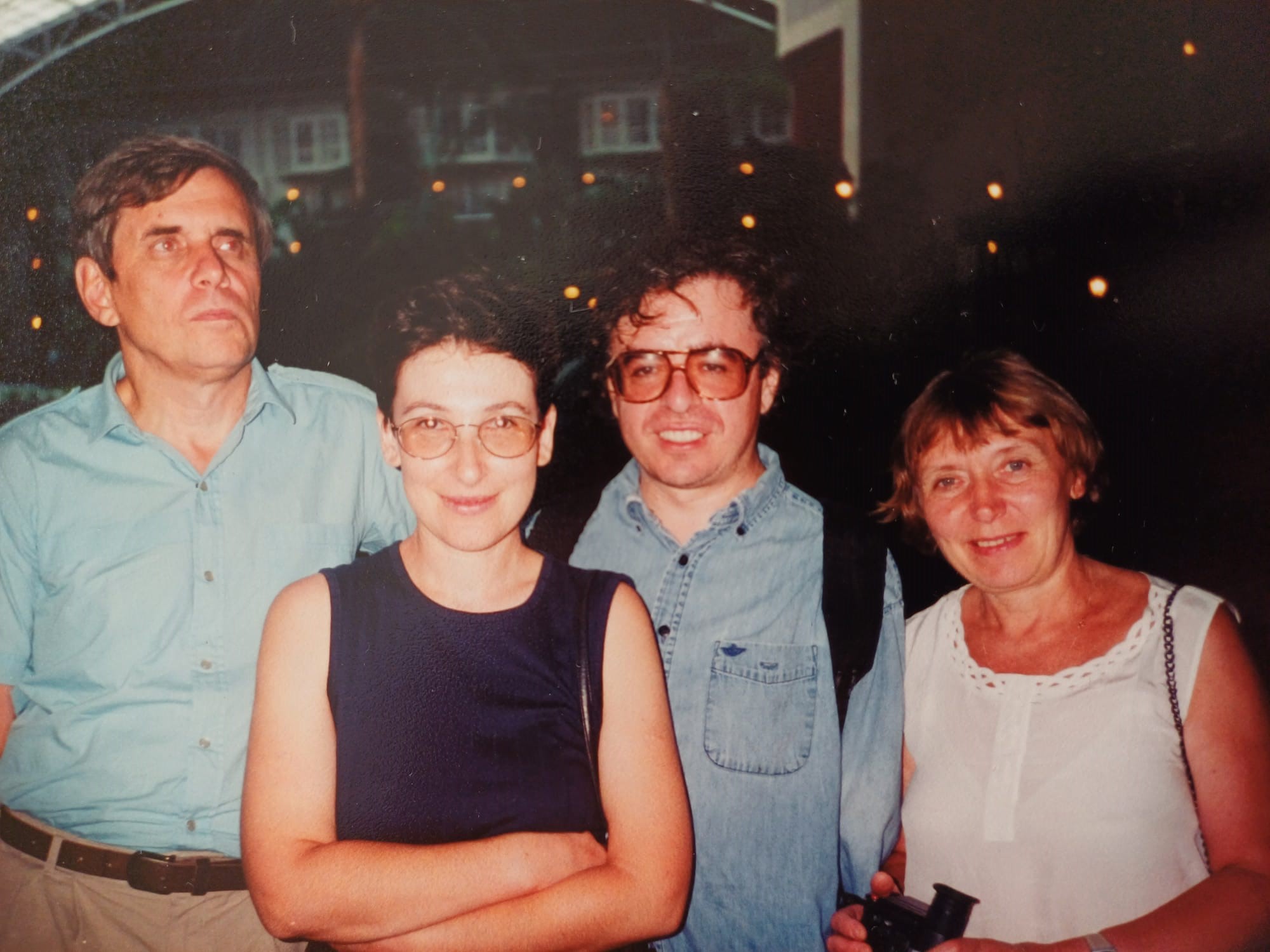}
		\caption{From left to right: Alexander Olshanskii, Olga Sapir, Mark Sapir and Tatyana Olshanskii in Nashville TN, 2000.}
		\label{fig:Olsh}
	\end{figure}

	Sapir gave an invited talk at the International Congress of Mathematicians in Madrid in 2006. He gave an AMS Invited Address at the American Mathematical Society Sectional Meeting in Huntsville, Alabama in 2008. 	For the outstanding contribution to mathematics, Sapir was elected the Fellow of the American Mathematical Society, inaugural class of 2012. A mathematical conference in honor of Sapir’s 60th birthday took place at the University of Illinois at Urbana–Champaign in 2017.
	
	Sapir founded the Journal of Combinatorial Algebra, published by the European Mathematical Society, and has served as its founding editor-in-chief since 2016. He has also served on the editorial board of several other journals.
	
	Mark was a passionate and dedicated mathematician who enjoyed discussing any problem, from the newest and most complex challenges to high school Olympiad ones. Both in Lincoln and in Nashville, Mark established a popular weekend school for talented young children. With enthusiasm and drive, Mark organized dozens of mathematical conferences and taught students in summer schools.
	He advised 7 graduate students and mentored many others.

	Mark was one of the first to offer his help and support to  immigrants and newcomers. 
	Apart from mathematics Mark talked about everything, and he had a great sense of humor. 	With stoicism, Mark faced the challenges of his later years. His death was a profound loss to all who knew him.
	
	 In this memorial, we have chosen to highlight two of his most significant accomplishments: his solution of Burnside-type problems for semigroups and his construction of S-machines. A brief mathematical overview of these topics follows. Additionally, we invited several of his colleagues and friends to share their reflections. 
 Their contributions, blending both professional and personal memories, appear below.

\subsubsection*{Burnside-type problems for semigroups}

Recall that a (semi)group $S$ is \textit{periodic} if each element in $S$ has finite order. A  stronger condition is that of local finiteness: a (semi)group $S$ is called \textit{locally finite} if each finite subset of $S$ generates a finite sub(semi)group. For groups, the general Burnside problem posed by William Burnside in 1902 asked whether every periodic group is locally finite.

A periodic (semi)group $S$ is of \textit{finite {exponent}} if  there exists a positive integer $n$ that bounds the order of each element of $S$. In the same 1902 paper, Burnside asked whether every group of finite exponent is locally finite; this is known as the bounded Burnside problem.

In the 1930s, a related question, called the restricted Burnside problem by Wilhelm Magnus in 1950, was asked for groups. It makes perfect sense for semigroups so we state it in the general form: Is there a constant $c(m,n)$ depending only on $m$ and $n$ such that $|S| \leq c(m,n)$ for every finite (semi)group $S$ with $m$ generators and exponent at most $n$? 

The three Burnside problems were extremely influential in the development of combinatorial group theory in the 20th century. 
The general and the bounded Burnside problems were solved in the negative in the 1960s by Golod and Shafarevich and by Novikov and Adian, respectively. 
 The positive solution to the restricted Burnside problem brought the Fields Medal of 1994 to Efim Zelmanov.

For semigroups,  analogs of these classical Burnside problems have negative solutions even if they are restricted to the variety of all semigroups with zero in which the square of each element is zero (so that the order of each elements is at most 2)—in this variety, there exist both infinite and arbitrarily large finite semigroups with 3 generators \cite{MH44}. On the other hand, restrictions of Burnside problems to some semigroup varieties may admit positive solutions; say, such is any variety of commutative semigroups. Can one ‘separate the sheep from the goats’, that is, find out in which varieties all periodic semigroups are locally finite and which contain infinite finitely generated periodic semigroups? This is the problem that was addressed in \cite{Sap87}. 

The varieties studied in \cite{Sap87} are of finite axiomatic rank, i.e.,  defined by identities  involving only finitely many variables. 
Sapir proved that, remarkably, for such a variety $V$, all periodic semigroups in $V$ are locally finite if and only if all  semigroups of finite exponent in $V$ are locally finite. Thus, in such varieties, the general and bounded Burnside problems   coincide. For non-periodic varieties $V$ (those containing an infinite cyclic semigroup),
 Sapir provided a complete solution, characterizing the varieties in which the general (or equivalently, the bounded) Burnside problem has a positive solution. For periodic varieties $V$, Sapir reduced the problem 
  to its version for the class of all groups in $V$. 
  
  Sapir's solution hinges on the sequence of Zimin words, \( \{Z_n\}_{n=1,2,\ldots} \), defined inductively by \( Z_1 = x_1 \) and \( Z_{n+1} = Z_n x_{n+1} Z_n \) in variables \( x_1, x_2, \ldots, x_n, \ldots \). Sapir’s results demonstrate that for a semigroup variety \( V \) defined by identities involving at most \( n \) variables, if \( V \) contains a nonperiodic semigroup or all groups in \( V \) are locally finite, then all periodic semigroups in \( V \) are locally finite if and only if 
   an identity  \( Z_{n+1} = W \), for some word $W$ distinct from $Z_{n+1}$, holds in $V$. 
  
  Note that the condition imposed on the variety $V$ in Sapir's solution  is very mild. Most ‘natural’ varieties
  are finitely based (i.e., defined by finitely many identities) and thus fall within the scope of Sapir’s results. For finitely based varieties, it is possible to algorithmically verify whether they satisfy an identity \( Z_{n+1} = W \) for some word $W$ distinct from $Z_{n+1}$. This provides an algorithmic way to determine if the general Burnside problem has a positive solution within a non-periodic finitely based variety.

Sapir's proof is a clever mixture of combinatorial arguments with those from the structure theory of semigroups. The original reasoning in \cite{Sap87} relied on a compactness argument from symbolic dynamics; another, fully constructive proof was published by Sapir in \cite{Sap91}. The latter allowed him to obtain  a solution for the restricted 
Burnside problem for finitely based semigroup varieties.
Indeed, Sapir gave an algorithm that, given any finite set $\Sigma$ of semigroup identities,  verifies whether or not the variety defined by $\Sigma$ has the ‘restricted Burnside’ property, that is, admits a uniform upper bound on the size of its finite semigroups with any fixed number of generators.

Reaching the current state of development of Burnside problems for groups required the efforts of many brilliant mathematicians over the span of 90 years. Remarkably, Mark Sapir brought Burnside problems for semigroups to a comparable level in just two papers written within 5 years.

\subsection*{S-machines}

	Let $G\cong\langle A; R\rangle$ be a finitely presented group. A word $u$ in the alphabet $A^{\pm 1}$ represents the identity in $G$ if and only if there exists a {derivation} $u=w_0\to w_1\to\dots\to w_d=1$, where every transition $w_i\to w_{i+1}$ involves either the cancellation or insertion of a pair of mutually inverse letters, or the deletion or insertion
of a subword  $r^{\pm 1}$ for $r\in R$. The\textit{ Dehn function} $d:{\mathbf N}\to{\mathbf N}$ of the group $G$ is defined to be the smallest function such that, for 
every word $w$ of length $\leq n$ representing $1$ in $G$,
there is a derivation $w\to\dots\to 1$ of length at most $d(n)$.

It is well known (and easy to see) that the word problem of a finitely presented group is solvable if and only if its Dehn function is recursive.
 It was observed in 
 \cite{SBR02} that if $d(n)$ is a recursive Dehn function of a finitely presented group $G$ then there is a ``trivial'' (nondeterministic) Turing machine 
 that checks whether a word $w$ is equal to $1$ in $G$,  by exploring all derivations of length at most $d(n)$ starting from $w$. This machine has nondeterministic time complexity  $O(d(n))$.

Nonetheless, there exist groups $G$ with a large Dehn function, for which the time complexity of solving the word problem is significantly lower. In these cases, the word problem can be solved efficiently by an optimized machine, while the ``trivial" machine may perform quite slowly.

The main result of \cite{BORS02} is that the word problem of a finitely generated group is decidable in polynomial time if and only if this group can be embedded (with linear distortion) into a group
with polynomial Dehn function. Moreover, every group with  word
problem solvable nondeterministically in time $T(n)$ (assumed to be superadditive) can be embedded with linear distortion into a group with Dehn function $n^2T(n^2)^4$.

Thus, if the word problem in a group can be solved rapidly by a sophisticated machine, this group can always be embedded in a finitely presented group where the word problem can be solved nearly as quickly by the ``trivial'' machine. 

This result also implies the existence of finitely presented groups with NP-complete word problem. As part of the proof, Birget, Rips and Sapir proved that, up to some mild restrictions,  Dehn functions above $n^4$ are the same as the time complexities of nondeterministic Turing machines \cite{SBR02}. In particular, they constructed  finitely presented groups with Dehn functions equivalent to $[n^{\alpha}]$ for ``almost arbitrary"
real $\alpha\ge 4$. (When \cite{SBR02} was in preparation, only integer
exponents $\alpha$ were known for Dehn functions of groups, the first irrational exponent was constructed in \cite{BB00}.)

The main ingredient in the long and intricate proof is Sapir's S-machines. Recall that there are various classical constructions for simulating a (deterministic) Turing machine by a finite  group presentation (Novikov, Boone, Britten, Higman, Aanderaa, ...), or for proving the Higman embedding theorem. However, all these constructions involve an exponential overhead; more specifically, the Dehn function of the finitely presented group is at least exponentially larger than the time complexity of the Turing machine (see \cite{MO85}). To obtain groups with polynomially bounded Dehn functions, Mark Sapir invented new, amazing machines called S-machines.

Unlike traditional Turing machines, whose rules depend on the current state of the machine and the letter(s) currently observed by the reading head, the nondeterministic S-machines are “blind”— their rules do  not depend on tape(s) content, although the head ``senses" when it reaches the end(s) of tape(s).

Sapir provided an intricate construction that  
for any given classical Turing machine 
$T$ with time complexity $T(n)$, produces a polynomially equivalent S-machine $S$. This means that $S$ and $T$
 recognize the same set of words, and that the time complexity of $S$ is bounded above by a polynomial function of $T(n)$. From this, a finite group presentation that simulates the S-machine $S$ and has  Dehn function polynomially-equivalent to $T(n)$ can be constructed. 
 
Even the creation of the initial phase of $S$, denoted $S_3$ in \cite{SBR02}, required considerable ingenuity; this component  can just  “distinguish zero from non-zero,” or, in other words, differentiate two letters in the tape alphabet. Sapir’s complete construction of an S-machine $S$ that is polynomially equivalent to $T$ spans several dozen pages, and to date, no one has managed to simplify it.

Sapir's S-machines have since become essential tools for solving several major open problems. 
Using these machines, Olshanskii and Sapir   introduced the first finitely presented counterexample to von Neumann's question on amenability: a finitely presented non-amenable group  with no non-cyclic free subgroups \cite{OS02}. 
They later
addressed Collins' 1976 problem by proving that a finitely generated group has a decidable conjugacy problem if and only if it embeds in a finitely presented group with the same property \cite{OS04}. In \cite{OS20}, 
they answered E. Rips' 1994 question by constructing a finitely presented group with  quadratic Dehn function and an undecidable conjugacy problem. This is significant since sub-quadratic Dehn functions imply hyperbolicity and thus decidable conjugacy.

	\section*{Victor Guba\footnote{Sadly, Victor Guba passed away on December 7th 2023.}}

	I remember the first time I saw Mark.  It was at the group theory seminar at Moscow State University, where  I was a graduate student at the time. Usually the talks were given by someone from our university, but sometimes we had outside speakers.  This talk was very unusual: most speakers would talk about their results in some narrow area, making it difficult for listeners to follow. But here, everything was clear and accessible, which surprised me.  After the talk, I asked our group leader: ``Who was that?'' - ``That was Sapir, from Sverdlovsk,''
	Olshanskii 
	 answered.

	Back then I did not know that Mark was a student of L.N. Shevrin, who, in addition to high scientific merits, was distinguished in his ability to give especially interesting talks, a skill that he passed on to his students.

	We were first introduced personally at the algebra conference in memory of A.I. Maltsev in Novosibirsk.  This was in 1989. It was a large meeting, with many famous participants from different countries.  When we met with Mark, I offered to tell him about Olshanskii's proof that the Burnside groups are infinite. I was fascinated by this topic at the time, and I remembered the entire paper well.  Mark was a very attentive and thoughtful listener. I told him about this outside the classroom, orally.  We wandered around the university campus for a long time, and I was able to lay out the proof in its entirety.  Later, many years later, the content of my story formed the basis of one of the chapters of Mark's book devoted to combinatorial algebra.

	After the conference, we were all sitting at the airport before leaving; Olga Sapir was there as well.  We had a lively discussion about some questions in semigroup theory, and I became very interested in this subject, which was then new to me.

	Then came the 90s.  The country underwent big changes.  Many mathematicians, including Mark, went to work abroad. In 1993 we met at a semigroup conference in York, England.  I remember an interesting conversation in which John Meakin also participated, where I first learned about diagram groups.  They were recently invented by Mark and John in one of their papers.  This was when I first heard about an open question about a diagram group over the simplest semigroup presentation with one relation, x=xx.

	And the following year, Mark invited me to visit him in the US. There, I lived for two weeks at his house.  At the time, Mark worked in Lincoln, Nebraska, along with Meakin and Margolis.  I really liked this small and quiet city, and in the evenings we discussed different math questions.  Once, after a conversation with Vesna Kilibarda, a student of Meakin and Sapir, I carried out the necessary computations and found that the very group we were talking about in York was the famous group of Richard Thompson. Now this is a well-known fact, but back then we were very surprised.  I immediately told Mark the full proof, and he was quite impressed.  It became clear that this class of groups deserved special study.  And this is what we did next.

	I cannot help but mention the hospitable and relaxed atmosphere in the family of Mark, Olga, and Jenya Sapir.  They were a friendly and welcoming family, and they made you feel right at home.  As hosts, they always took care of you, and everything you needed.

	In 1995, at the invitation of the university, I came to Lincoln to teach for one semester.  That's when we started writing up our work on diagram groups.  We had already obtained some of these results, but in the course of our work, we made more and more new discoveries.  Here I must mention one of Mark's significant qualities: he was incredibly inventive in posing new questions.  This is something that does not come naturally to me: I like to explore that which is already present.  But here, while we were writing, as a result of our constant discussions, so many new topics and chapters appeared.  And so, we wrote a book that was over 100 pages long, which was soon published by the AMS.  I must say that there are now several hundred references to it. It was written in a very short period of time, and was the result of an enthusiastic and fruitful collaboration.

	Mark's passion and hard work were absolutely infectious to those around him. Collaboration with him was always effortless, despite the fact that the amount of work to be done was enormous.

	After this, we corresponded regularly, sending each other several long e-mails a day.  All this correspondence has been preserved; it contains quite a lot of valuable information.  Sometimes I re-read individual letters and wonder at how everything happened.  I think we complemented each other perfectly in many respects, so the work went very efficiently. Many more articles were written about diagram groups in the process of this intensive communication.

	Here it is impossible not to recall one joke that came from colleagues.  In the field of diagram groups, we have been called the ``founders'', using the same word in Russian that is used to refer to Marx and Engels. This could not help but amuse the people who lived in the USSR and remembered those ideological times.  Needless to say that Mark had great respect for Karl Marx and his legacy.  I had a different way of looking at things, and as a result, we had various debates between us.  Mark always tried to express his point of view with arguments, with links to facts or texts.  As I now understand, he believed in the possibility of reorganizing society on a foundation of reason, within some systematic doctrine.  It was funny that in Russian the words MarkS and Marx sound the same, and the first of these was often chosen by Mark as a nickname on various websites.

	In the summer of the same year, 1995, Mark visited Russia for the first time since leaving for the United States.  At first he lived for some time at my house in Vologda, after which he visited his hometown (then already called Yekaterinburg), and then we went together to St. Petersburg for a semigroup conference in honor of E.S. Lyapin.  It was a perfectly organized conference, in which Slava Kublanovsky played a significant role. He was already working in the private sector by then, but at the same time did not stop doing mathematics.  The combination of these two activities is rare in and of itself.  But he also proved some outstanding results, which Mark spoke highly of, and on the basis of which he defended his doctoral dissertation.

	Since then, Mark and I met many more times at conferences in various parts of the world: in the US, Portugal, France, England.  Mark's talks at all these places were distinguished by their originality and their ability to draw in his listeners, to get them to understand the most important points.  Everything was presented quickly, energetically, accessibly, and at just the right level of rigor.  He always gave a live performance, without using prepared slides, as is now common.  It was noticeable that a lot of our colleagues appreciated his talks, including those working on other topics.  I must say that Mark had a very broad mathematical outlook in general, as well as extensive knowledge of many different areas of algebra.  It is no wonder he became the author of so many articles and monographs that go far beyond  any narrow, specialized field.

	What should be especially noted is Mark's  cheerfulness and his sense of humor.  Everyone who knows him remembers his characteristic rolling laugh, with which he liked to react to so many situations.  I also remember a particular story: one of his American students complained that he rarely hears words of approval from Mark. To me this is very understandable, since I know Mark's character, and I understand that it was not typical for him to express any ``obligatory'' official words.  And he reacted to this in his usual manner, placing an audio file on his home page for those who wish to listen to ``words of encouragement'' from him.  There, set to music, in a monotonous, robotic voice, are words like ``good job!'', ``great job!'', ``I like your homework!'' and the like.

	After this, Mark moved from Lincoln to Nashville, to work at a larger university, and I was invited several times to teach there.  All those years our joint work continued.  But, strange as it may seem, our work progressed most intensively when we went into in absentia mode, when we were in different countries.  However, personal meetings and discussions were also very meaningful.  Among other things, we spent time with each others' families: Mark first had his son, Yasha, who was a very lively and active boy, and then his youngest daughter Rachel - our children also spent time with one another.  And in Nashville he soon began to work with A.Yu. Olshanskii, with whom Mark had developed a close collaboration.  Together they wrote many significant and voluminous articles in the field of combinatorial group theory.  In general, I must say that Mark was able to simultaneously work on a large number of new projects, which from the outside could not help but impress.  Moreover, many hardworking people do this sort of thing as a job, viewing it as a series of tasks that one must be particularly disciplined and diligent to accomplish.  But to Mark, it was just a part of his life, and seemed to take no special effort, despite the colossal volume of work that he produced.  What was started was always completed. Before my eyes, he wrote a long seminal work together with Rips and Birget (I was even its reviewer).  It was there that Mark introduced the concept of their famous S-machines, which became the basis of many new and important things that came after.

	After some time, Mark underwent a complex brain surgery.  Preparing for it, Mark, in his usual cheerful manner, said something like the following to his colleagues: ``and if I pass, then in such and such a year you will hold a conference in my memory''. Thank God, at the time everything went well, and for many years after this (more than 10 years) Mark did not experience serious health problems.

	The last time we met in person was in 2016 in Denmark, at a conference in Odense.  In those years, I already traveled abroad quite rarely.  It was nice to see the familiar faces of colleagues after a long break.  Mark was still full of energy and cheerfulness.  In fact it seemed that over the years he did not change at all - this also applied to his everyday habits and demeanor.

	A few years later, I learned from Olshanskii that Mark had serious new health problems.  He did not travel any more, but continued to actively engage in mathematics.  We corresponded at that time, although not so often. At the same time, he managed to follow the new results of his colleagues, and devoted a lot of time to his journal, where he was the editor in chief.  Back then, somehow, I did not want to think about sad things, and there was hope that Mark would overcome his illness this time too.  However, this did not happen, and in October 2022, the sad news reached all of us...

	Remembering Mark, I always think about the bright and joyful moments of his life, and how he knew how to captivate those around him with his enthusiasm.  Objectively, now you understand what an outstanding person he was in so many respects. This, despite the fact that all this was taken for granted when one interacted with him, and also despite the fact that he was an extremely modest person.  But one way or another, he continues to live on in our memory, in the  things we do that are connected to him, in his children.  May his soul be eternal.
	
	\begin{figure}[h]
		\centering
		\includegraphics[width=0.47\textwidth]{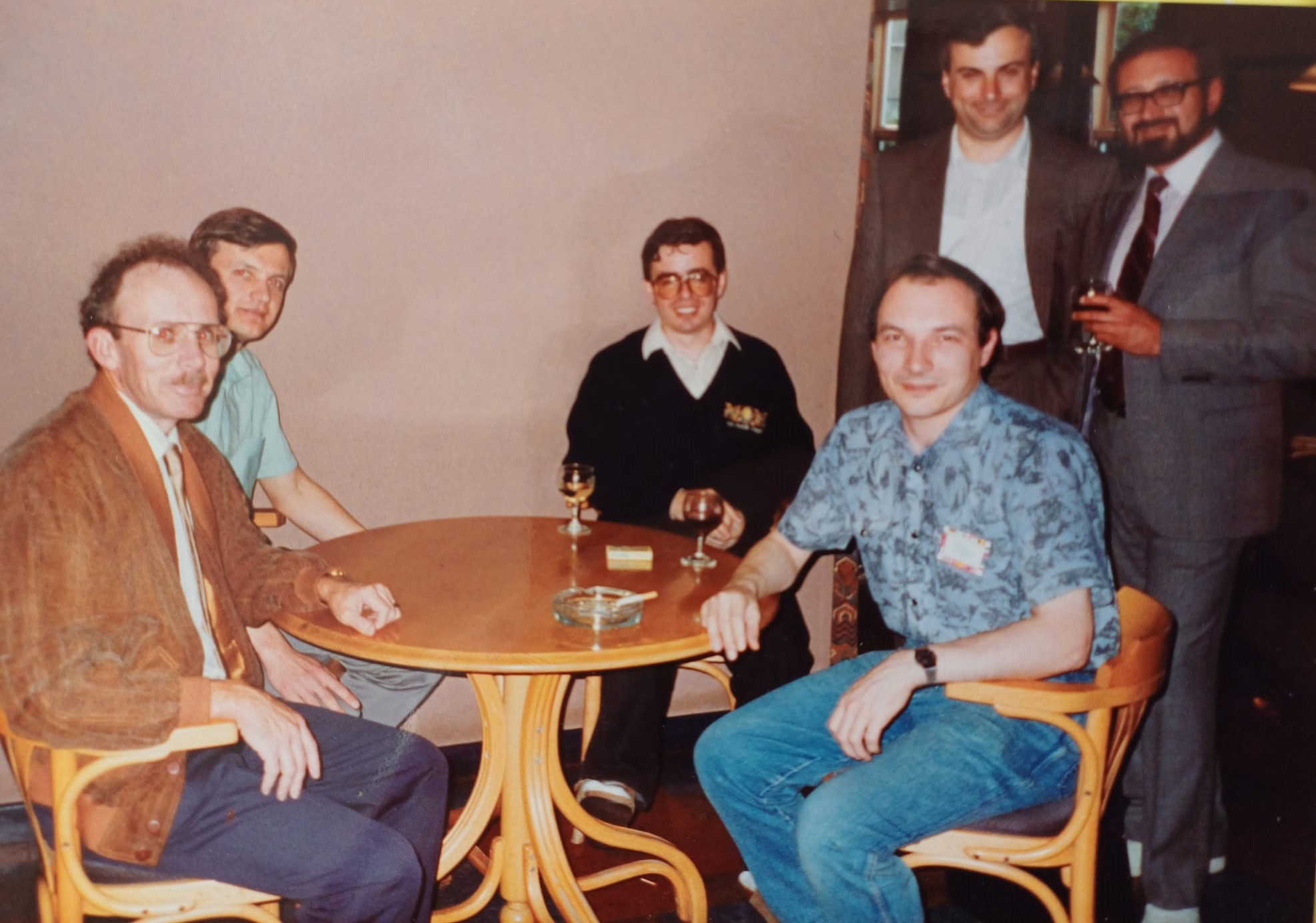}
		\caption{From left to right: (sitting) John Meakin, Andrei Kelarev, Mark Sapir and Victor Guba, (standing) Misha Volkov and Simon Goberstein, in NATO ASI, York, England 1993.}
		\label{fig:con}
	\end{figure}
	
	\section*{Olga Kharlampovich}
	
	Those who knew Mark Sapir are greatly saddened by his untimely death. Mark was an old friend from our undergraduate student years at the Ural State University in Sverdlovsk (now Yekaterinburg) Russia. Mark was a few years ahead of me but we shared the same  supervisor Professor Lev Shevrin. Mark had a remarkable enthusiasm for math and devoted himself to math education. As a PhD student he was involved in teaching math to undergraduates as well as high school students. He actively participated in the ``winter school" organized by the math department, where students immersed themselves in mathematics, solved problems, enjoyed skiing, singing, and celebrating together.
	I vividly recall the festive atmosphere of the ``winter school" that we as well as the other participants enjoyed so much!  
	
	While we were in Sverdlovsk we were active in the seminar  ``Algebraic Systems" headed by Prof. Shevrin, the seminar attracted around 30 participants from various institutes (later called universities) in Sverdlovsk. In the period 80-91 Mark focused on algorithmic and other problems in varieties of semigroups and associative algebras and gave numerous talks in the seminar. 
	He obtained a complete description of finitely based non-periodic semigroup varieties where relatively finitely presented semigroups have decidable word problem. Every relatively finitely presented semigroup in such a variety is representable by matrices over a field. Mark also described varieties of semigroups where every finitely generated semigroup is residually finite and obtained many other classification results for semigroup varieties. I worked on similar questions for groups.
	
	In 1989 our joint work on the survey ``Algorithmic problems in varieties" started. Our two families became very friendly; our children were of the same age. However, our paths then diverged.  In 1990, I moved to McGill University in Montreal, Canada, while Mark accepted a position at the University of Nebraska, Lincoln, in 1991. He collaborated with John Meakin at the University of Nebraska. They developed in a 1993 paper the notion of a diagram group, based on finite semigroup presentations. Mark further developed this notion in subsequent joint papers with Victor Guba.   Important examples of diagram groups are Thompson groups. 
	
	We met again in 1992 in Lincoln, shared our first impressions of North American life and restarted the work on the survey. We were mainly concerned with varieties of classical algebras, like groups, semigroups, associative and Lie algebras.  But we  aimed to present results in the most
	general form so when possible we formulated statements for arbitrary universal
	algebras. In addition to algorithmic problems we were dealing with their
	``neighbors and relatives" like residual finiteness, the Higman embedding property,
	the finite basis property etc. We also discussed the computational complexity of
	solvable algorithmic problems. Our survey was published in 1995 in the International Journal of Algebra and Computation.  Mark's contributions to the journal did not end there; he later became its managing editor and held the position until 2015.
	
	Mark joined Vanderbilt University as a professor of mathematics in 1997. At that point his interests shifted towards geometric group theory. Here, he achieved some of his most renowned and significant results, notably in two papers published in the Annals of Mathematics in 2002. Collaborating with Jean-Camille Birget and Eliyahu Rips in the first paper and with Birget, Rips, and Alexander 
	Olshanskii in the second, Mark provided a complete description of all growth types of Dehn functions of finitely presented groups. Additionally, they demonstrated that a finitely presented group has the word problem solvable in NP if and only if it embeds into a finitely presented group with a polynomial Dehn function. 
	These are truly remarkable papers. 
	
	His interests extended further, as he and Cornelia Dru\textcommabelow{t}u developed the asymptotic cone approach to study relatively hyperbolic groups. Mark also authored several exciting papers with Alexander Olshanskii, further enriching the field of geometric group theory.
	
	Mark Sapir's untimely departure leaves a void in the world of mathematics and in the hearts of those who admired him. His legacy lives on in his groundbreaking works and the profound impact he had on his colleagues and students. 
	
	\section*{Simon M. Goberstein}
	
	In the Spring of 1976, as a graduate student at the Saratov State University, I attended a conference on semigroups at the Leningrad Pedagogical Institute. One of the participants was Émil Golubov, an Associate Professor at the Ural State University whom I had met earlier. I saw Émil in the evening before the first day of the conference, and he told me that he had come to Leningrad with his co-author, a brilliant third-year undergraduate student, Mark Sapir, who in Émil's own words was a rising star with a great future in mathematics. That was the first time that I heard the name Mark Sapir, and I met Mark the following day when he gave a talk on his joint work with Golubov on residually finite semigroups. The results presented by that 19-year-old undergraduate were so impressive that it was clear that he is indeed a rising star! Moreover, Émil's prediction that this young student will have a great future in mathematics was spot-on. The breadth and depth of Mark Sapir's knowledge and understanding of mathematics were truly remarkable, and he has made seminal contributions to several areas of algebra and related fields, especially, to semigroups, to geometric group theory, and to combinatorial algebra. Mark also had a rare ability of presenting the most difficult topics in a very accessible form, and he was extremely generous in sharing his ideas with colleagues and students. It is impossible to give a comprehensive description even of his major results in a short article. Hopefully, some of his co-authors will highlight the most significant aspects of their joint work with Mark. Here I would like to confine myself to the recollection of my interactions with Mark and his family. Thus this note is more personal than professional.      

	For a few years after our first meeting, Mark and I saw each other only occasionally. I emigrated from the Soviet Union in January of 1980 with my wife and our 6-year-old son, and in 1981, we settled in Chico, California. In 1989, I saw Mark again at a conference in Berkeley, and he asked me if I could find a position for him at my department. It took me some time to make all the appropriate arrangements, and at another conference at Monash University in July of 1990, I informed Mark that the invitation to him and his wife, Olga, will be sent in a couple of months. Finally, in January of 1991, Mark, Olga, and their little daughter, Jenya, arrived in Chico, and Mark got a position of a Visiting Associate Professor at Chico State. That was the beginning of our friendship that lasted until Mark's untimely death on October 8, 2022. 

	In the Spring semester of 1991, Mark was assigned to teach two sections of Finite Math for Business. The students in such courses usually have a negative attitude toward mathematics. However, Mark quickly adjusted to the situation and was able to build a great rapport with his students. He was also quite interested in computer science education and was actively promoting the special software package, Roo and Robby, that he had developed back in Russia. Mark was an avid cyclist, and we found a bike for him shortly after their arrival in Chico. He was very impressed by the peaceful and quiet life in a small university town. He even wrote to one of his friends that there is no crime in Chico and it was safe to leave his bike unlocked near the entrance to their apartment. Unfortunately, he was too optimistic about that -- the morning after he had mailed the letter, his bike was stolen. Of course, our department faculty were eager to help and donated two bikes to Mark and Olga (one with a child seat for Jenya). I also taught Mark to drive a car, using my old Chevette (which had an automatic transmission but no power steering). He learned very quickly, passed the driving test with flying colors, and got his first driver license. 

	Of course, Mark and I talked a lot about mathematics and started working together on some problems concerning the lattice of quasivarieties of inverse semigroups. Mark was also continuing to revise a comprehensive survey, {\it Algorithmic problems in varieties} (written jointly with O. Kharlampovich), which was published shortly thereafter. It was a fortunate coincidence that Professor Alexander Olshanskii had a visiting position at UC, Berkeley in the Spring semester of 1991. We invited him to come to Chico to give a colloquium talk and then again to visit us for a few weeks. Mark Sapir and Alexander Olshanskii knew each other back in Russia (Professor Olshanskii was an `official opponent' at the defense of Mark Sapir's Ph.D. thesis), and they were both happy to see each other again. It is quite likely that their fruitful cooperation started in Chico during Olshanskii's visit. Of course, it continued with great success and intensity when they both became colleagues at Vanderbilt.

	Finally, I must add that our families were also very close. We had a number of joint trips, visited each other many times, especially, on holidays and for family celebrations. Mark had an excellent sense of humor, and we all had great time together. He and I enjoyed talking about literature, history, and current events, and it was especially gratifying since our tastes and views were very similar. In the last months of his life, Mark and I often talked on the phone, and I know that his spirit was strong to the very end. Mark Sapir's passing is a great loss to the mathematical community and a deep personal loss for me.

	\section*{Stuart W. Margolis}

	I first met Mark in 1987 at a conference in Szeged, Hungary. This was the period of Glasnost and it was the first time that semigroup theorists in the West were able to meet their colleagues in the East in such numbers. It was very exciting for all of us. Given the volume and depth of Mark’s work that I knew about, I had assumed that he was an older mathematician. I was very surprised when the young man who entered the lecture hall was identified as Mark. 
	
	We were able to communicate via email over the next year. The next time we saw each other was at a conference on semigroups at the University of California-Berkeley organized by John Rhodes in 1989. At that time, Mark was interested in immigrating to the United States. John Meakin and I worked very hard to convince the University of Nebraska-Lincoln to give Mark an appointment in the mathematics department in 1991, where he remained until moving to Vanderbilt University in 1997. This was certainly one of the greatest hires made at UNL. With Mark, John Meakin, Jean-Camille Birget, myself and our students we had a very fruitful and exciting semigroup seminar at UNL during those years. Among the people who visited us for extended times were Victor Guba, Eliyahu Rips and Alexander Olshanskii. This interaction led to the two famous papers with Mark, Birget, Rips and the second also with Olshanskii on isoperimetric functions, isodiametric functions and computational complexity \cite{SBR02},\cite{BORS02}. This work was the first steps that Mark made into Geometric Group Theory where he did most of his groundbreaking research in the 21st century. 
	
		\begin{figure}[h]
		\centering
		\includegraphics[width=0.47\textwidth]{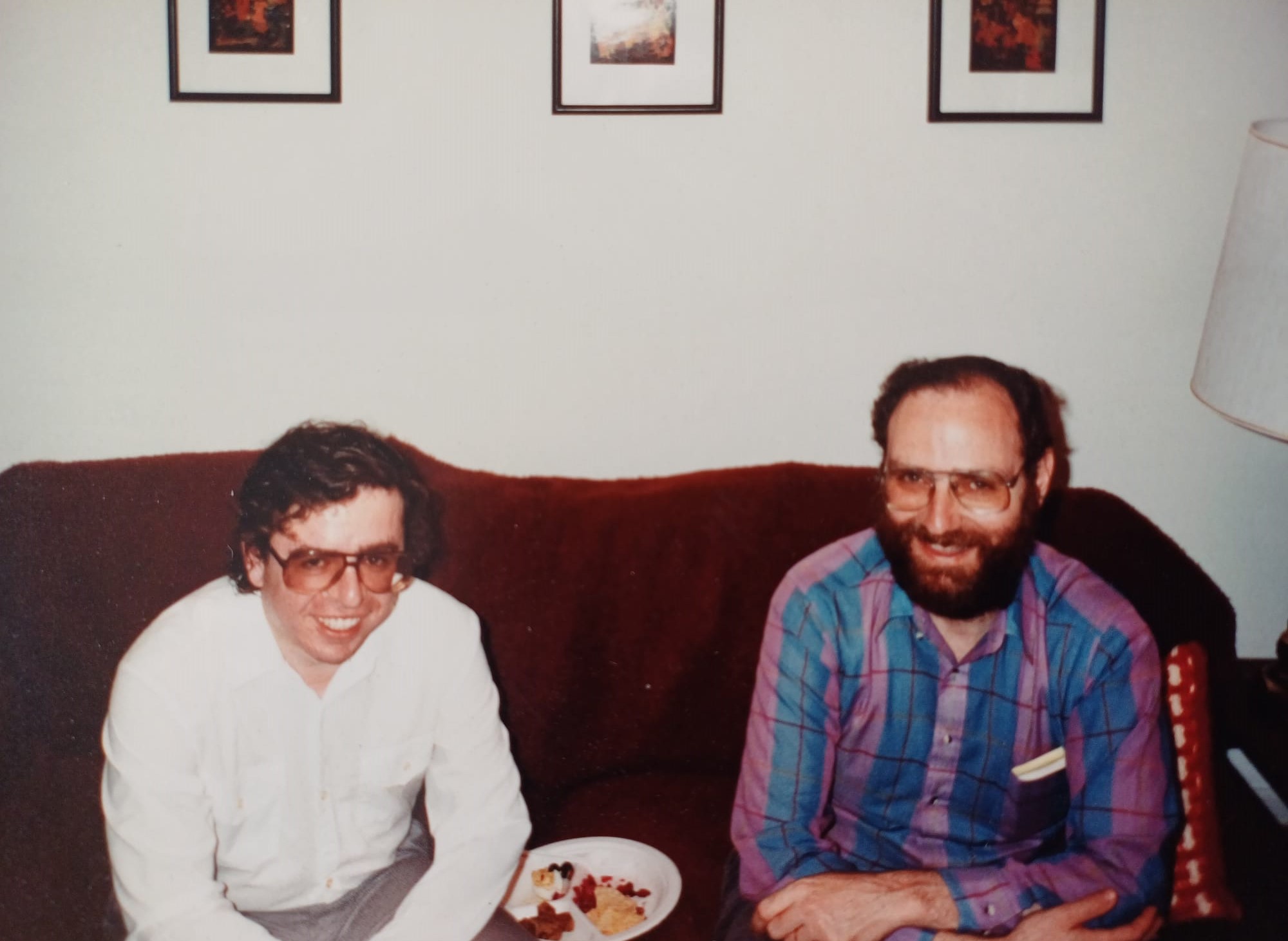}
		\caption{Mark Sapir and Stuart Margolis in Lincoln, Nebraska, May 1993.}
		\label{fig:mar}
	\end{figure}
	
	As is well-known, Mark was a student of Lev Shevrin, the founder of the famous school of semigroups at the Ural State University. Mark's thesis and career up to 1997 was mainly in various aspects of semigroup theory. Not surprisingly, Mark made major contributions in this field as well. It was an honor for me to be a co-author of six papers with Mark, along with other collaborators. To give a flavor of the kind of research that we were involved in at the University of Nebraska in the 1990s, I'll summarize a few. These results have been significantly developed since in the last 30 years.
	
	The paper \cite{MSW01} studied extension problems for finite automata of partial bijections. The most important example of such automata are the Stallings automata associated to finitely generated subgroups of free groups \cite{Sta83}. 
	The question studied is if one can decide if a collection of partial bijections on a finite set $Q$ can be extended to a finite permutation group in a pseudovariety $V$ of finite groups on some set $P$ containing $Q$. This turns out to be deeply related to the pro-$V$ topology on the free group.  The paper 
	\cite{MS95} shows that every finite semigroup belongs to a locally finite finitely based quasi-variety of semigroups. That is, there are no finite inherently non-finitely based semigroups with respect to quasi-identities. This is in sharp contrast with the case of varieties and identities where many finite semigroups are inherently non-finitely based with respect to identities. The proof uses methods of symbolic dynamics developed by Mark going back to his Ph.D. thesis.
	
	One of my favorite results that appeared in the paper \cite{HKMST97}  showed that the problem of deciding if a finite semigroup is a subsemigroup of a finite 0-simple semigroup is undecidable. Finite $0$-simple semigroups are one of the basic building blocks of all finite semigroups. They were determined by Rees in the 1940s and have a quite transparent structure as matrix semigroups over a group. Thus it was very surprising to find out that one can not decide if a semigroup embeds into a finite 0-simple semigroup. Moreover, the problem is intimately related to the positive first order theory of finite semigroups. Once again, the proof was dependent on tools developed by Mark earlier in his career, this time what he called split systems. Along with John Meakin, we wrote a long survey article \cite{MMS95} on the state of the art on algorithmic problems in groups, semigroups and inverse semigroups up to 1993. In 1991 when John Rhodes, John Meakin and I started the International Journal of Algebra and Computation (IJAC), we asked Mark to be one of the founding Editors. Later he became the Managing Editor of the journal. Mark had very high standards and this helped the journal attract papers of the highest quality. 
	
	Mark was very generous in sharing his ideas and his methods with me and other collaborators. It was a lot of fun to work with him. He had a wonderful sense of humor. The years 1991-1996 in which we were colleagues at the University of Nebraska-Lincoln were among the most fruitful and exciting of my career in a large part due to Mark. It was wonderful to have Mark and Olga, of course a powerful mathematician herself and his young family, including his young daughter Jenya who went on to become a mathematician, adapt to Lincoln, Nebraska, which must have seemed to be something like the land of Oz compared to their native Yekaterinburg (then Sverdlovsk). Mark was one of the greatest mathematicians that I have personally known and he enriched my life both personally and mathematically. May he be remembered for blessings.
	
	\section*{Lev Shneerson}
	
	I knew Mark for more than 40 years. We met for the first time at the end of the 70's - my friend Evgeniy Sukhanov told me a year before about a brilliant young student actively participating in the internationally famous Professor L. N. Shevrin's Seminar in Sverdlovsk (now Yekaterinburg).
	
	We met with Mark many times during the Algebraic Conferences and we had common points of interest in Combinatorial Semigroup Theory.
	
	Mark's talks were always exceptional. He combined a very deep knowledge of various areas of Modern Algebra and other branches of Mathematics with phenomenal intuition. He was able to obtain bright results using elegant proofs with links from one area to another. I remember how during Mark's talk, one of the leading Hungarian algebraists, László Márki, who sat next to me exclaimed with admiration that it was worth coming to the Conference just to listen to this talk.
	
	Mark's solutions of the classical Burnside type problems in important classes of semigroup varieties obtained in the 80s have amazing and unexpected applications not only for semigroups. For example, he simultaneously established that the set of semigroup identities of a nilpotent group $G$ is finitely based (in the class of all semigroups) if and only if $G$ is abelian or has a finite exponent. This proposition became a very important complement to the well-known results of A.I. Malcev and  B.H. Neuman-T.Taylor who proved independently that for any positive integer $c$, the variety of all nilpotent groups of class nilpotency less than or equal to $c$ can be defined by a certain semigroup identity.
	
	I often remember a joint trip with Mark and Misha Volkov to Australia to the Conference in honor of G.B. Preston in Melbourne and the Workshop on Semigroups in Sydney in 1990. Mark suggested to look at the preliminary text of my first talk. He came to my room late at night and carefully read everything - giving very helpful advice regarding slides that I never used before.
	
	We had very interesting discussions after my talk at the Conference in York in 1993 during which he posed a conjecture about connections between the growth and unavoidable words, and during the Conference in St. Petersburg in honor of E.S.  Lyapin in 1995. 
	
	Mark was one of the first who offered his help and supported me when I immigrated to the US. I will always remember him.

	\section*{John Meakin}
	
	I first met Mark Sapir at an international conference on semigroups in Szeged, Hungary in 1987. This was one of  a series of important conferences in Szeged that enabled interaction between mathematicians from the west with their counterparts from the Soviet Union and other eastern European countries. I believe that this was Mark's first international conference outside of Russia. He was one of  a small group of brilliant mathematicians from Lev Shevrin's ``Algebraic systems" seminar in Sverdlovsk (now Yekaterinburg) who traveled to Szeged 
	 for this conference. Even at that early stage of his career, it was obvious that Mark was rapidly developing into a leading figure in the field, with powerful and surprising results in the general area of varieties and quasivarieties of semigroups.
	
	Another opportunity to interact with Mark presented itself two years later when we met at a workshop on monoids at Berkeley in August 1989. Mark obtained last minute permission to participate in the workshop and he arrived at Berkeley with an open date return airline ticket on Aeroflot but no return reservation, so he needed to spend several additional days in Berkeley. Arrangements were made for him to stay with mathematicians in the Berkeley area until he could secure a return airline reservation, and so we were able to engage in wide ranging conversations about problems of mutual interest. In fact our first joint paper, linking congruences on free monoids with certain inverse submonoids of polycyclic monoids, had its genesis in our discussions at the time of the Berkeley workshop. Thus began a personal friendship and mathematical collaboration that extended for many years.
	
	I was fortunate to be able to arrange a faculty position for Mark at Nebraska in 1991.  Mark moved to Lincoln with his wife Olga and their young daughter Jenya (now on the faculty in mathematics at Binghamton) in August 1991. The family lived in Lincoln until 1997 when Mark accepted a position at Vanderbilt.  Their son Yasha was born in Lincoln. The family made the transition to life in the USA: Olga completed her PhD in mathematics, Jenya began her education in the Lincoln Public School system and Mark learned to drive a car (an interesting experience that he attacked with his characteristic intensity and sense of humor, even learning to fix the inevitable occasional dents and scrapes to his car himself). There were warm family dinners, gatherings with mathematical visitors and graduate students, a vibrant seminar and free-ranging conversations (about mathematics, politics, western popular culture, films, any number of topics).
	
	Mark had a major impact on mathematical life at Nebraska, not only through the power of his brilliant research, but also through his commitment to teaching and student learning at all levels. While at Nebraska, he organized a Saturday mathematics school for 7th through 10th grade school students and a popular series of Saturday morning schools for young programmers (grades 4 through 7), based on an educational software package that he had helped develop before he left Russia. He also created a WebBook on linear algebra for undergraduates. His lectures, to students and to colleagues in seminars and at professional conferences, were always exceptional models of clarity and insightful exposition. Mark  established strong collaborative research relationships with colleagues in Lincoln (particularly Jean-Camille Birget, Stuart Margolis and myself), with a focus on strengthening  connections between geometric group theory and semigroup theory. He attracted a series of prominent invited seminar speakers in both areas and he helped organize one of a series of international research conferences at Nebraska with this focus.
	
	\begin{figure}[h]
		\centering
		\includegraphics[width=0.47\textwidth]{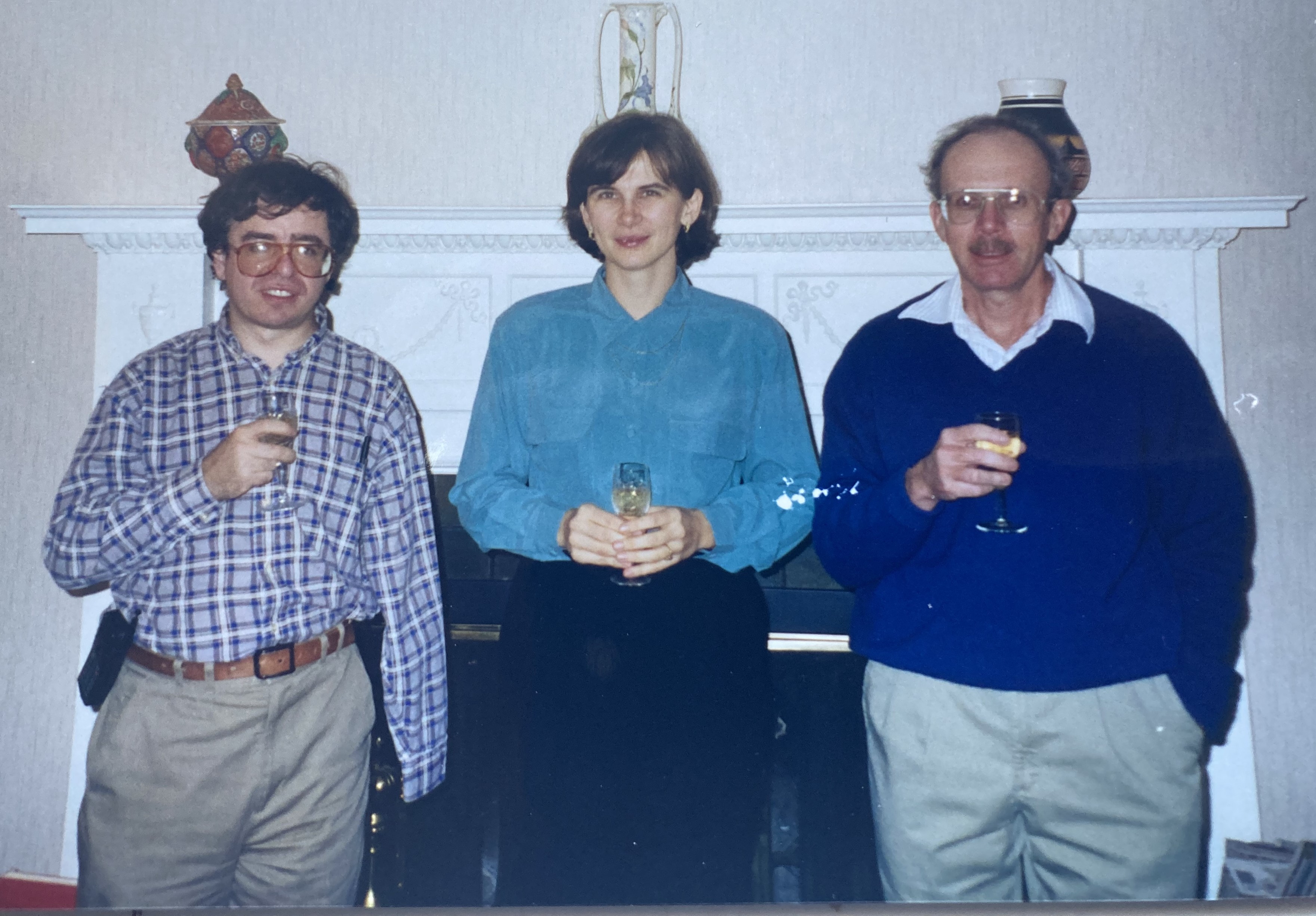}
		\caption{Mark Sapir, Tanya Jajcayova and John Meakin in Lincoln, Nebraska, December 1997.}
		\label{fig:Lincoln}
	\end{figure}

	I remember Mark as a fearless mathematician, willing (in fact eager) to discuss any mathematical problem with anyone, intellectually generous in his advice to students and colleagues but demanding in his expectation of commitment to the discipline. I remember his drive, his ingenuity, his optimism and good humor in the face of many challenges, his passion for mathematics, his inspiring intellect, the depth of his thinking. It was a great privilege to have had him as a colleague and to have had the opportunity to work with him. He was a mathematician of the first order, and he will be remembered for his deep contributions to algebra, particularly in semigroup theory, geometric group theory and combinatorial algebra.

	\section*{Ilya Kapovich}
	
	I first met Mark Sapir in January 1994, when I was a graduate student in New York City. Mark gave a talk in the New York Group Theory Seminar at the CUNY Graduate Center. The seminar was run by my PhD advisor Gilbert Baumslag, and I greatly benefitted from attending talks there by many outstanding mathematicians from around the world. Mark Sapir was one of them. At the time Mark was still mainly working in semigroup theory, but he was already beginning to switch his interests to infinite groups. Mark was also an invited speaker at a geometric group theory summer school in Banff  in August 1996. That's where I really came to first appreciate Mark's remarkable mathematical insights. He was able to explain many deep results in group theory, such as the existence of a finitely presented group with an unsolvable word problem, in an amazingly intuitive and accessible way, with geometric pictures that really spoke to me and made the subject much easier to understand. I was also struck by Mark's idea to move beyond the Turing machines and to use more flexible computational devices, such as Minsky machines, when exploring and generalizing various classic results from group theory.  This approach found spectacular applications in the joint work of Mark with Olshanskii, Rips and Birget, published in two seminal papers in ``Annals of Mathematics" in 2002, on relating the computational complexity of a finitely generated group $G$ and the isoperimetric functions of finitely presented groups into which $G$ can be embedded as a subgroup, and on characterizing the possible types of functions that can occur as isoperimetric and isodiametric functions of finitely presented groups.  
	
	I followed Mark's work closely ever since then and I was fortunate enough to hear many of his talks and mini-courses in subsequent years, and to read his beautifully written papers. Although Mark came from a rather algebraic and combinatorial mathematical background, he fully embraced and mastered a wide variety of tools from other areas of mathematics and applied these tools to geometric group theory, including ideas from topology, geometry, dynamics, analysis, and so on.   For example, Mark's work with Cornelia Dru\textcommabelow{t}u on the use of asymptotic cones and $\mathbb R$-trees considerably clarified and advanced the theory of relatively hyperbolic groups and their generalizations. Mark championed a truly holistic approach to mathematics, often with great success. One of my favorite papers of Mark is his 2005 article in ``Inventiones", joint with Borisov, where they use unexpected and rather nontrivial connections with algebraic geometry to prove that the mapping torus of an injective endomorphism of a finitely generated linear group $G$ is residually finite.  Previously this result was unknown even for the case where $G$ is a finite rank free group.
	
	From his time back in the Soviet Union, Mark retained a genuine passion for mathematical outreach and for getting children excited about learning mathematics. He had been a long term faculty participant of the Canada/USA Mathcamp, an annual immersive summer program for talented middle and high school students.  During the Spring 2017 semester, when Mark was a G.A. Miller Visiting Professor at the University of Illinois at Urbana-Champaign, he showed me some of the problems and educational materials from his Mathcamp work, and Mark's enthusiasm for getting kids excited about mathematics was truly impressive.  
	
	Mark was also an inspirational and energizing force for the young researchers working in group theory. From my own time as a graduate student and a postdoc, I still remember that, whenever I met Mark, he was always genuinely interested in what I was working on and in what else beginning researchers in geometric group theory were doing. When Mark asked ``What's new?", he really meant it, and his intellectual curiosity was infectious.
	
	Mark's deep and broad work greatly enriched and significantly advanced group theory and adjacent areas of mathematics. But what I, and I am sure many others, will miss most about Mark is his energy, enthusiasm and love of mathematics.

	\section*{Eugene Plotkin}

	My memories of Marik as a mathematician are tightly intertwined with memories of our various get-togethers, conversations, and trips. All this conjures in my mind the image of a remarkable person and an outstanding scholar.
	
	Before I begin, it is important to mention that I will sometimes be calling him Mark and other times by his pet name Marik, which corresponds much more closely with the warm feeling I remember him with and with the image of him I have in my mind.
	
	Marik and I met in the early 80s at one of the all-Union algebra conferences. Before that, however, we must have run around the same streets, attended similar kindergartens, and played in similar playgrounds, as we were both born around the same time in the same city - Sverdlovsk, now Yekaterinburg, in the heart of the Urals.
	
	At the conference, it was impossible not to be impressed by him. He immediately grasped the essence of any mathematical discussion, showed a complete understanding of whichever algebraic problem was being discussed. As for Marik's own studies, he was at the time working on semigroups, as you would expect from a student of Lev Shevrin, but one couldn't help but notice the breadth and depth of his mathematical potential.
	
	When we became friends and began to discuss various mathematical topics, it became very clear that in addition to being naturally gifted, he was also extremely hard-working and possessed an ability to focus on a problem completely. And he never cut himself any slack.
	
	One day, many years later, we were walking along the beach of the Pacific Ocean near San Diego. The weather was great, and I began to make tourist plans.
	
	``Marik," I said, ``Shall we go to one of the national parks?''
	
	``No,'' said Mark, with unexpected resolve. ``I can't.''
	
	``But why? Look how beautiful it is, and I think we have time!''
	
	``I can't,'' Mark said again. ``Alexander Yurievich (Olshanskii) 
	 said not to.''
	
	``Stop fooling around!'' I was starting to get annoyed. ``What do you mean, he said not to?''
	
	``I can't,'' he repeated. ``I need to work, or Olshanskii will be upset, he's waiting.''
	
	They were indeed working on a paper together and Mark simply couldn't—and didn't want to—divide his attention between that and something else. But at the same time, he didn't want to hurt my feelings by turning down my plan just like that, so instead he blamed it all on Olshanskii, who was, in fact, the nicest, gentlest person. It eventually became clear to me that trying to push on him was no use. He was a person of integrity. It would have been easier to make a river flow backward than to persuade Mark to change his mind. Over the years, we argued a lot about various mathematical and near-mathematical matters and Mark wouldn't give an inch defending his personal view of the subject.
	
	Once I asked Mark to review a short paper submitted for a mathematical volume I was editing. To say the review was scathing would be an understatement. However, I had my own thoughts on the matter, so I convinced Mark to ask the author to rewrite the paper. In the end, Mark approved revised version number 49! And he was right.
	
	In fact, all our friends were aware of Marik's very particular, often extreme views on certain mathematical results. I once tried, in a friendly chat, to convince him that some result was absolutely outstanding and deserved high praise. Marik turned abruptly and said: ``There is no such result! It needs verification. It needs to be sent to a peer-reviewed publication. Let them find reviewers, let them get it rewritten and set it out in a better way. Until then—there is no such result, and I don't want to even talk about it!''
	
	He was steaming, it was not easy to calm him down. He was right, of course. However, I did not like the sharp way in which he was putting it and was starting to get wound up myself: ``So what are you suggesting? Why don't you propose something positive instead?''
	
	Maybe it was then, or maybe earlier, that Mark had the idea of having long and complex mathematical texts checked and reviewed by a group of scholars. And it worked! It was indescribably hard to find three or four mathematicians who would agree to dissect and review a work together. But Mark was able to find convincing arguments and overcome objective obstacles. In other words, he could be persuasive, especially if the matter was worthy. I can think of at least two truly remarkable works that were under this ruthless revision process. It did not always result in publication success, but it would always result in much higher quality texts.
	
	Now I would like to step aside from mathematics-related memories. Marik and I shared two other passions, two topics we could discuss forever: mushrooms and table tennis. Marik loved ``mushroom-hunting'' like nothing else. I received photos of mushrooms from him from all over the world to identify. Once he sent me something, my immediate reaction to which was: ``Whatever you do, don't eat it!'' Marik was well aware: ``Don't worry,'' he wrote, ``this is purely scientific interest.''
	
	Once Marik paid us a visit at our home in Tel Aviv early in the morning. It was December and it was raining heavily and incessantly. My wife Tanya served us a breakfast of freshly picked slippery jacks. There was a sparkle in Marik's eye.
	
	``Where is this from?''
	
	``From the forest, obviously.''
	
	``Can we go?''
	
	``Marik, you have a colloquium in Jerusalem at four.''
	
	``That's enough time,'' Marik said confidently.
	
	I thought about it. There probably was enough time.
	
	``Tanya,'' I said. ``Could you please find something waterproof for Marik.''
	
	Tanya looked at the two of us dubiously: the difference in size was all too apparent. In the end she produced some trousers of mine and a rain jacket and having put these on, with a mushroom basket in his hand, Marik began to look like a theorem from iterated small cancellations theory, where an example of a group is built, unimaginable in our ordinary flat world. We went to my special place in Ben Shemen. The mushrooms were abundant, so was the rain; an hour later we were completely soaked, but our rain-washed mushrooms sat cozily in our baskets.
	
	``Mark, let's go,'' I said. ``It's all clay around, the road will get washed out, we'll be late for your talk at the Hebrew University.''
	
	Mark could barely be seen underneath the huge rain jacket, all I could hear was his happy giggles. We just about made it for his talk. Mark walked into the colloquium room on the first floor just as he was, in my huge trousers and rain jacket and with a basket full of mushrooms.
	
	``I was in the forest! I picked a ton of mushrooms! I made it in time for my talk and now I'm about to show you results just as beautiful as these mushrooms.''
	
	He came up to the whiteboard, took off my old wet clothes and without any thought or hesitation got straight to it. He wasn’t giving a talk, he was flying!
	
	Another non-mathematical passion we had in common was table tennis. On one of his visits to Israel we set up a date to go and play. I found a table in our friends' house, brought rackets—and off we went. Mark was playing well, sometimes really, really well. I had to try quite hard and was forced to remember all I had been taught in tennis school a long time ago. Finally, I started to win and, probably unconsciously, stopped serving professional spinning serves in order not to overload his weak hand. Mark noticed it and shouted: ``Play normally!'' Things like that are impossible to forget. He never made allowances for himself. Not in tennis, not in mathematics. And many times, he won!
	
	Mark had many friends, he liked good company, and he was good at maintaining friendships. He liked humor, he appreciated a good joke, he even had a subscription to the Israeli humor magazine ``Beseder'' and loved to quote fresh jokes. But at the center of his interests was always mathematics. I had the impression that he always had all questions and answers to hand, that he never had any difficulty remembering old results, because all of that had been processed and kept neatly inside. Once we were discussing construction of a finitely presented Burnside group. This problem had interested Mark since he was young, and he considered it one of the most difficult problems in group theory. However, it has never been solved and remains open to this day. In the end, Mark wrote in his 2007 overview that out of all finitely presented ``monster'' analogues he believed most of all in the existence of a finitely presented infinite torsion group. I think he may have been right: Mark's intuition was almost never wrong. We will just have to wait a bit longer...
	
	I would like to finish this little memoir by saying how grateful I am for having met Mark and for our friendship.
	
	\begin{figure}[h]
		\centering
		\includegraphics[width=0.48\textwidth]{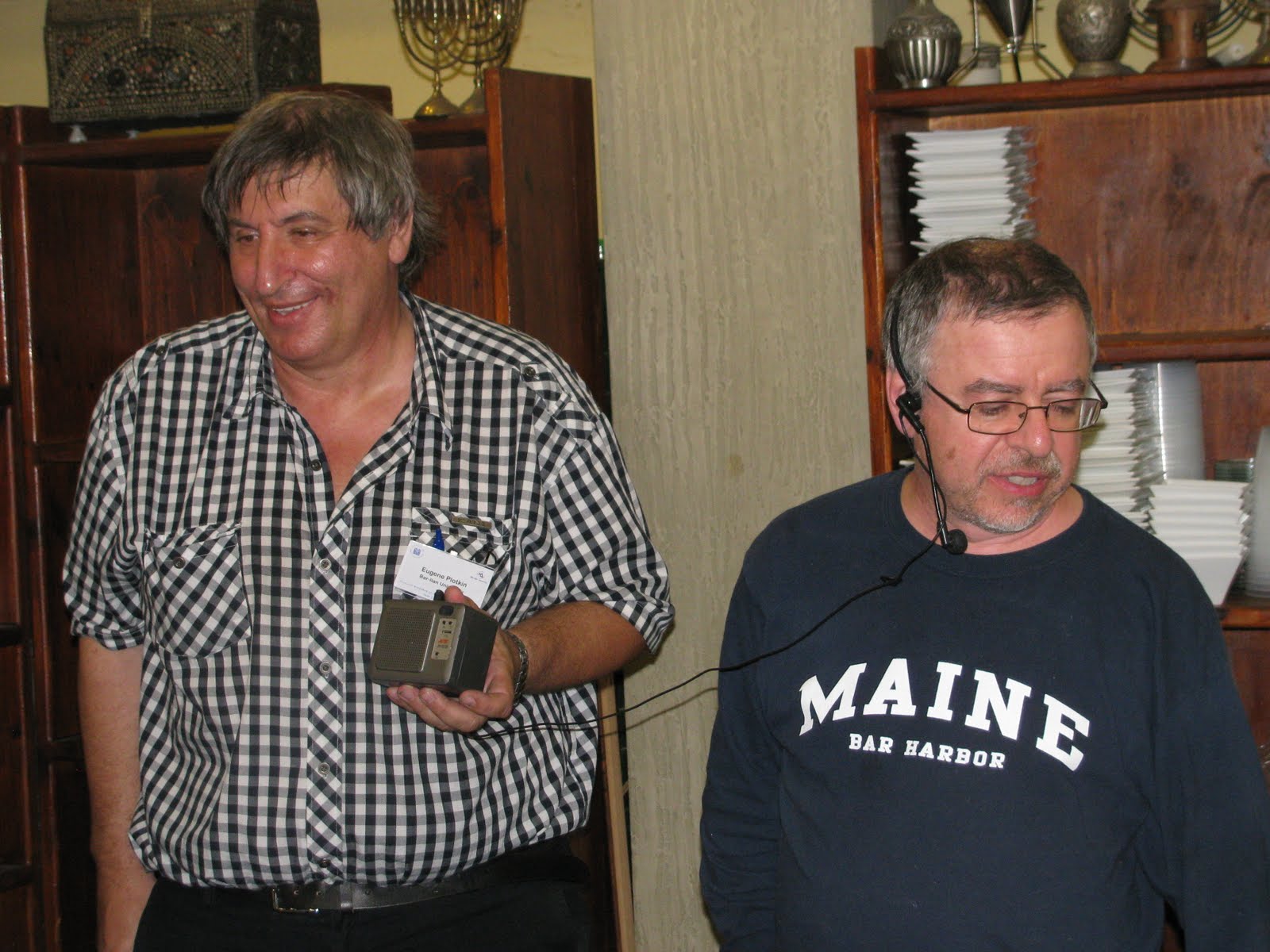}
		\caption{Eugene Plotkin and Mark Sapir in Ramat-Gan, Israel, 2013.}
		\label{fig:Plotkin}
	\end{figure}

	\section*{Michael Mihalik}
	
	Mark Sapir’s mathematical legacy is sealed with deep theorems and innovative techniques that are well documented by others. I am tasked with writing a more personalized memorial. Mark was a positive force in any seminar,  colloquium or conference that he attended. As a speaker, Mark had a gift of  delivery that made the most complex situation seem understandable. As a  member of the audience he always seemed to have a deep observation or  penetrating question – often accompanied by his classic chuckle – that  injected an element of excitement and energy into the room. For a speaker, it  was a special opportunity to be on the other end of one of Mark’s comments.  One of the things I will most miss about Mark is hearing his standard  comment following an adequate response to one of his questions: ``I  understand”.
	
	Outside of mathematics, Mark was a dedicated father and husband.
	 Hiking and bicycling were among his favorite physical activities. Mark loved the taste of smoked wild duck and goose and very much wanted to try hunting. I was hesitant at first, because it’s a bit dangerous and conditions can be brutal. You must be on the lake in the winter, an hour before sunrise. The boat ride is cold, damp and you simply cannot make a mistake. So why do it? Sunrise on the lake during full migration is spectacular. The pastel color of the sky reflecting on the water and the symphony of migrating birds is a unique experience. When Mark took his first goose, I may have been as excited as he was, but probably not. Sharing that with a dog and a friend is a rare type of experience and I am fortunate to have shared it with Mark. He didn’t hide his satisfaction on or off the lake. It was clear that he really did ``understand”. I miss him dearly and I am sure I will always remember his companionship when I’m on the lake or in a seminar without him.
	
	\section*{Samuel Corson}
	
	Mark Sapir was my PhD supervisor.  My first memory of him is at a conference in the late 2000s; the first real acquaintance was a couple years later at Vanderbilt University.  Those who meet Mark in a professional setting quickly recognize the speed and breadth of his mind.  He happily attended colloquia on various topics which are distant from his specialty, at which he would ask questions of impressive insight.  When he attended a talk on algebra his mind would very quickly see the essential ideas of a proof, sometimes divulging them before the speaker wished to do so.   Mark’s contributions to the corpus of mathematics are well-known.  Though highly-regarded as a group theorist, many forget that his interest in group theory began well into his career.  His insights into algebraic combinatorics and algorithmic problems lead to deep results in, for example, the complexity of the word problem and the structure of asymptotic cones.
	
	Mark’s talks would often give the impression that a solution was not difficult.  He had a clever way of revealing the essentials of an argument through an appropriate analogy.  Before delivering the punchline of a proof, the side of his mouth would curl into a smile and he would say, “It is easy.”  His speed and facility of presentation were so impressive, in fact, that three-quarters of the way through a semester-long introductory algebra course he finished the standard set of lecture notes.
	
	In personal interactions he was a kind and straightforward man, quick to chuckle at a joke.  His sense of humor is transparent to those who visit his website, especially his “Words of Encouragement”.  During downtime in his office he would play online chess.  As a supervisor he was encouraging and ambitious; he gave me papers to referee which I would not entrust to myself at that point in my career.  I feel honored to have known him.
	
	\section*{Arman Darbinyan}
	
	Mark Sapir was one of those professors at Vanderbilt University whose impact has had a significant effect on the trajectory of my research and mathematical career. I learned about him as one of the deepest and most original minds in geometric group theory from the very first days of my PhD career, but it was during my second year at Vanderbilt that I got a chance to interact with him more closely. My systematic meetings with Mark started when I took a reading course with him, during which we were introduced to combinatorial algebra based on his, at that time not yet published, book of the same name. This course helped me to appreciate the beauty and intuition of the subject that later became my research area. Mark had a unique ability to explain the intuition and zest behind seemingly impenetrable and technically involved results in very accessible ways. Through interacting with him I learned about the creative process that shapes intuition and then transforms it into theorems and theories. Most importantly, I shaped a perception of mathematical aesthetics and learned more about the magic of the interplay between syntax and semantics in mathematical theories. Later, I took several other reading courses with Mark and had many occasions to enjoy the journey into new mathematical areas under his guidance. A theorem coauthored by him and proved using his theory of S-machines, that bridges geometric group theory and the theory of computational complexity, through its depth and beauty guided me in my research that resulted in my dissertation and several consequent papers.

	Besides being a brilliant mathematician and mentor, Mark was also an admirable human, whose very subtle jokes, culture, erudition in literature and poetry, discussions of chess, and friendly support were invaluable stimuli during my PhD years and beyond. Through his passion, energy, charisma, and initiative, Mark was also playing a significant role in community building among geometric group theorists. I am very grateful to Mark for making me feel welcome in that community.
	

	
	\section*{Credits}
	Figure \ref{fig:Mark} is courtesy of John Russell, Vanderbilt University. 
	
	\noindent Figures \ref{fig:Olsh},\ref{fig:con} and \ref{fig:mar} are courtesy of Olga Sapir. 
	
	\noindent Figure \ref{fig:Lincoln} is courtesy of John Meakin. 
	
	\noindent Figure \ref{fig:Plotkin} is courtesy of Stuart Margolis.
	
	\vspace{1.77em} 
	\noindent Jean-Camille Birget,\\
	Department of Computer Science, Rutgers University - Camden, 
	Camden, New Jersey 08102, USA.\\	
	e-mail: birget@camden.rutgers.edu. \\
	
	\noindent Gili Golan,\\
	Department of Mathematics, Ben Gurion University of the Negev, 
	Beer Sheva 8410501, Israel.\\
	e-mail: golangi@bgu.ac.il\\
	
	\noindent Alexander Olshanskii, \\
	Department of Mathematics, Vanderbilt University, Nashville, TN 37240, USA. \\
	e-mail: alexander.olshanskiy@vanderbilt.edu\\
	
	\noindent Mikhail Volkov,\\
	Institute of Natural Sciences and Mathematics, Ural Federal University 620000 Ekaterinburg , Russia.\\
	e-mail: m.v.volkov@urfu.ru\\
	
	\noindent Victor Guba, \\
	Vologda State University, 15 Lenin Street, Vologda 160600, Russia.\\
	e-mail address: gubavs@vogu35.ru\\
	
	\noindent  Olga Kharlampovich,\\
	Department of Mathematics and Statistics, CUNY Hunter College and Graduate Center, 695 Park Avenue, 10065, New York, USA.\\
	e-mail: okharlam@hunter.cuny.edu\\
	
	\noindent Simon M. Goberstein,\\
	Department of Mathematics and Statistics, California State University, Chico, CA, 95929-0525, USA\\
	e-mail: sgoberstein@csuchico.edu\\

	\noindent Stuart W. Margolis,\\
	Department of Mathematics, Bar-Ilan University, Ramat Gan, 5290002, Israel\\
	e-mail: margolis@math.biu.ac.il\\

	\noindent Lev Shneerson,\\
	Department of Mathematics and Statistics, Hunter   College of CUNY, 695 Park Avenue, 10065, New York, NY, USA.\\
	e-mail: lev.shneerson@hunter.cuny.edu\\

	\noindent John Meakin,\\
	Department of Mathematics, University of Nebraska, Lincoln, Nebraska 68588, USA.\\
	e-mail: jmeakin@unl.edu\\

	\noindent Ilya Kapovich,\\
	Department of Mathematics and Statistics, Hunter College of CUNY, 695 Park Avenue, New York, NY, 10065, USA.\\
	e-mail: ik535@hunter.cuny.edu\\

	\noindent Eugene Plotkin,\\
	Department of Mathematics, Bar-Ilan University, Ramat Gan, 5290002, Israel\\
	email: plotkin@math.biu.ac.il.\\

	\noindent Michael Mihalik,\\
	Department of Mathematics, Vanderbilt University, Nashville, TN 37240, USA. \\
	e-mail: michael.l.mihalik@vanderbilt.edu\\
		
	\noindent Samuel Corson,\\
	Matematika Saila, UPV/EHU, Sarriena S/N, 48940, Leioa - Bizkaia, Spain\\
	e-mail: sammyc973@gmail.com\\
	
	\noindent Arman Darbinyan,\\
	School of Mathematical Sciences, 
	University of Southampton,
	University Road 
	Southampton 
	SO17 1BJ 
	United Kingdom.\\
	e-mail: a.darbinyan@soton.ac.uk\\

\end{document}